\documentclass[12pt]{article}
\usepackage{amsmath, amsthm, amssymb, stmaryrd, fullpage}
\usepackage{amsrefs}
\usepackage[all]{xy}

\numberwithin{equation}{subsection}
\newtheorem{theorem}[equation]{Theorem}
\newtheorem{lemma}[equation]{Lemma}
\newtheorem{cor}[equation]{Corollary}
\newtheorem{prop}[equation]{Proposition}
\theoremstyle{definition}
\newtheorem{defn}[equation]{Definition}
\newtheorem{hypothesis}[equation]{Hypothesis}
\newtheorem{example}[equation]{Example}
\newtheorem{remark}[equation]{Remark}
\newtheorem{convention}[equation]{Convention}

\def\AAA{\mathbb{A}}

\def\FF{\mathbb{F}}
\def\GG{\mathbb{G}}

\def\QQ{\mathbb{Q}}
\def\RR{\mathbb{R}}
\def\ZZ{\mathbb{Z}}
\newcommand{\calC}{\mathcal{C}}

\newcommand{\calO}{\mathcal{O}}
\newcommand{\calS}{\mathcal{S}}
\newcommand{\gotho}{\mathfrak{o}}

\newcommand{\bv}{\mathbf{v}}
\newcommand{\bw}{\mathbf{w}}

\def\dual{\vee}

\DeclareMathOperator{\alg}{alg}
\DeclareMathOperator{\Aut}{Aut}

\DeclareMathOperator{\coker}{coker}

\DeclareMathOperator{\Hom}{Hom}

\DeclareMathOperator{\Mor}{Mor}
\DeclareMathOperator{\perf}{perf}

\DeclareMathOperator{\Spec}{Spec}

\begin{document}

\title{On the geometry of $p$-typical covers in characteristic $p$}
\author{Kiran S. Kedlaya \\ Department of Mathematics \\ Massachusetts
Institute of Technology \\ 77 Massachusetts Avenue \\
Cambridge, MA 02139 \\
\texttt{kedlaya@math.mit.edu}}
\date{January 16, 2006}

\maketitle

\begin{abstract}
For $p$ a prime, a $p$-typical cover of a connected scheme on which $p=0$ is a finite
\'etale cover whose monodromy group (i.e., the Galois group of its
normal closure) is a $p$-group.
The geometry of such covers exhibits some unexpectedly pleasant
behaviors; building on work of Katz, we demonstrate some of these.
These include a criterion for when a morphism induces an isomorphism of
the $p$-typical quotients of the \'etale fundamental groups,
and a decomposition theorem for $p$-typical covers of polynomial rings
over an algebraically closed field.
\end{abstract}

\section{Introduction}

Let $p$ be a prime number.
A finite \'etale cover of a connected scheme on which $p=0$ is 
\emph{$p$-typical} if the monodromy group of the cover (which for a connected cover
coincides with
the Galois group of the normal closure) is a $p$-group.
The geometry of such covers exhibits some unexpectedly pleasant
behaviors; the purpose of this paper is to briefly expose a few of these.
This is in part to dispel the notion that one can only ever prove
meaningful results about the tame
(prime-to-$p$) quotient of the \'etale fundamental group.

For instance, Katz has shown \cite{katz-canonical}*{Proposition~1.4.2}
that if $R$ is a connected ring in which $p=0$, then the categories of
$p$-typical covers over $R[t^{-1}]$ and over $R((t))$
are equivalent, via the evident base change functor.
In other words, if $\pi_1^p$ denotes the maximal pro-$p$ quotient of
the \'etale fundamental group $\pi_1$ (where basepoints are suppressed
throughout this introduction for notational simplicity), 
then the natural homomorphism
$\pi_1^p(R((t))) \to \pi_1^p(R[t^{-1}])$ is a bijection.
We give a natural generalization of Katz's theorem
(Theorem~\ref{T:p-limit}), which characterizes more generally when
one connected affine scheme of characteristic $p$ looks like a
limit of a diagram of others from the point of view of constructing 
$\pi_1^p$.
Here is a sample result (Example~\ref{E:Laurent series}): 
if $k$ is an algebraically
closed field of characteristic $p>0$, then 
\[
\pi_1^p(k[t,t^{-1}]) \cong \pi_1^p(k[t]) \times \pi_1^p(k[t^{-1}]).
\]
(The analogous statement for $\pi_1$ is false:
the left side has nontrivial prime-to-$p$
quotients whereas the right side does not. Note also that in general,
neither $\pi_1$ nor $\pi_1^p$ commutes with products, so one cannot
replace the right side with a single fundamental group of $\AAA^2_k$.)

We also look more closely at $p$-typical covers of affine toric varieties,
including of course ordinary affine spaces. Our main results in this direction
(Theorem~\ref{T:discrete} and its corollaries, notably
Theorem~\ref{T:splits2}) assert that the $\pi_1^p$ of any affine toric
variety can be written as an inverse limit of $\pi_1^p$'s of one-dimensional
varieties, or even of affine lines. This requires the use of some
auxiliary ``height functions'' to measure the complexity of
$p$-typical covers; we can describe some simple examples of such functions,
but only \emph{a posteriori} (Theorem~\ref{T:height splits}).

\subsection*{Acknowledgments}

Thanks to Dan Abramovich for helping track down an error in a previous 
version of this paper, and to the referee for some helpful suggestions
concerning the exposition.
The author was supported by NSF grant DMS-0400727.

\section{$p$-typical covers}
\label{sec:artin-sch}

In this chapter, we introduce the notion of a $p$-typical cover
and prove a strong generalization of Katz's canonical extension
property for such covers (Theorem~\ref{T:p-limit}).
We first fix some notational conventions for the whole paper.
\setcounter{equation}{0}
\begin{convention}
Throughout this paper, fix a prime number $p$.
By a ``$p$-group'', we will mean a finite group whose order is a power of $p$.
Standard facts about $p$-groups, which we
will use without further comment, include the following.
\begin{enumerate}
\item[(a)] The center of any nontrivial $p$-group is nontrivial.
\item[(b)] Any maximal proper subgroup of a nontrivial $p$-group is normal of
index $p$.
\end{enumerate}
By a ``$p$-ring'', we will mean a ring in which $p=0$; likewise for
``$p$-field'' or ``$p$-domain''. Similarly, by a ``$p$-scheme'',
we will mean a scheme in whose ring of global sections one has $p=0$.
\end{convention}

\subsection{The \'etale fundamental group}

We first recall the notion of the \'etale fundamental group from
\cite{sga1}*{Expos\'e~V} (with some notation as in
\cite{katz-canonical}*{Section~1.2}).

\begin{convention} \label{conv:conn}
Throughout this section, let $X$ be a connected scheme, and
let $\overline{x}$ be a geometric point of $X$,
i.e., a morphism $\Spec k^{\alg} \to X$ in which $k^{\alg}$ is an algebraically
closed field.
\end{convention}

\begin{defn}
Let $\calC_X$ denote the category of finite \'etale covers
of $X$; note that $\calC_{\overline{x}}$ may be identified with the
category of finite sets.
Then the pullback functor $F_{\overline{x}}: \calC_X \to \calC_{\overline{x}}$
is represented by a pro-object $P$ of $\calC_X$. Let 
$\pi_1(X, \overline{x})$ denote the automorphism group of 
$F_{\overline{x}}$, i.e., the group of pro-automorphisms of $P$.
\end{defn}

\begin{remark}
Replacing $\overline{x}$ by another geometric point $\overline{y}$
does not change the abstract structure of the group $\pi_1(X)$.
However, there is no \emph{canonical} isomorphism
$\pi_1(X,\overline{x}) \to \pi_1(X,\overline{y})$; the choice of
such an isomorphism constitutes the choice of a ``chemin'' (``path'').
\end{remark}

\begin{defn}
Let $X$ be a connected
scheme, let $E \to X$ be a finite \'etale cover,
and let $\overline{x}$ be a geometric point of $X$.
Then the profinite group $\pi_1(X,\overline{x})$ acts continuously
on $E_{\overline{x}}$, and the image is well-defined up to group
isomorphism.
We call it the \emph{monodromy group} of $E$.
\end{defn}

\begin{defn}
If $E \to X$ is a connected finite \'etale cover, 
there is a unique minimal connected finite
Galois (\'etale) cover $E' \to X$ which factors through $E$; it is the maximal
cover fixed by the kernel of the map $\pi_1(X, \overline{x}) \to
\Aut(E_{\overline{x}})$. Consequently, the Galois group of this cover
is precisely the monodromy group of $E \to X$. This cover is called the
\emph{normal closure} (or \emph{Galois closure}) of $E \to X$; it coincides
with the usual field-theoretic definition when $X = \Spec k$.
\end{defn}

\subsection{$p$-typical covers}

We now extract the $p$-typical part of the fundamental group.
Throughout this section, we retain Convention~\ref{conv:conn}.

\begin{defn}
A \emph{$p$-typical cover} of $X$ is a finite \'etale
cover $E \to X$ whose monodromy group is a
$p$-group; if $S/R$ is a ring
extension whose corresponding cover $\Spec S \to \Spec R$ is $p$-typical,
we say $S$ is a \emph{$p$-typical extension} of $R$.
Note that the fibre product and the disjoint union of $p$-typical covers are
$p$-typical.
If $E$ is connected and $p$-typical over $X$, then $\deg(E\to X)$
is a power of $p$: namely, this degree is the index in the monodromy group
of the stabilizer of any geometric point of $E$. 
\end{defn}

\begin{lemma} \label{L:ptyptrans}
If $E \to X$ and $E' \to E$ are finite \'etale covers with $E$ connected, 
then $E'\to X$ is $p$-typical if and only if $E'\to E$ and $E\to X$ are both
$p$-typical.
\end{lemma}
\begin{proof}
Choose a geometric point $\overline{x}$ of $X$ and a geometric
point $\overline{y}$ of $E_{\overline{x}}$. Let $G$
be the monodromy group of $E'\to X$, identified with the image of
$\pi_1(X, \overline{x})$ in $\Aut(E'_{\overline{x}})$, and let
$H$ be the monodromy group of $E'\to E$, identified with the image of
$\pi_1(E, \overline{y})$ in $\Aut(E'_{\overline{y}})$.
Then $H$ is the stabilizer of $\overline{y}$ within $G$.

On one hand, if $E'\to E$ and $E\to X$ are $p$-typical, then $H$ is a
$p$-group, $G$ acts transitively on the geometric points of $E_{\overline{x}}$
(since $E$ is connected), and so $\#G = \#H \cdot \deg(E\to X)$ is a
$p$-power. Hence $E'\to X$ is $p$-typical.

On the other hand, if $E' \to X$ is $p$-typical, then 
$G$ is a $p$-group, as then must be $H$, so $E' \to E$ is $p$-typical.
Meanwhile, the monodromy group of $E \to X$ is a quotient of $G$, since
any element of $\pi_1(X, \overline{x})$ fixing $E'_{\overline{x}}$
must in particular fix $E_{\overline{x}}$. Hence $E \to X$ is also
$p$-typical.
\end{proof}

\begin{defn}
Let $\calC_X^{p}$ denote the subcategory of $\calC_X$ consisting of 
$p$-typical covers. Again, the fibre functor $F^p_{\overline{x}}: \calC_X^{p}
\to \calC_{\overline{x}}$ is represented by a pro-object of
$\calC_X^p$, whose group of pro-automorphisms
coincides with the automorphism group of $F^p_{\overline{x}}$.
We call this group $\pi_1^p(X,\overline{x})$ and refer to it as
the \emph{$p$-typical fundamental group} of $X$; the
inclusion $\calC_X^{p} \hookrightarrow \calC_X$ induces a surjection
$\pi_1(X, \overline{x}) \to \pi_1^p(X,\overline{x})$, under which
$\pi_1^p(X,\overline{x})$ is identified with the maximal
pro-$p$ quotient of $\pi_1(X,\overline{x})$.
\end{defn}

\subsection{$p$-typical covers and Artin-Schreier towers}

We will mainly be interested in $p$-typical covers of $p$-schemes;
these can be studied using Artin-Schreier towers. 

\begin{defn}
For $G$ a finite group (viewed as a constant group scheme over
$\Spec \ZZ$) and $X$ a scheme, a \emph{$G$-torsor} over $X$ is
a finite \'etale cover $E \to X$ equipped with an action of $G$, which
\'etale locally on $X$ is isomorphic to $X \times G$ (the 
\emph{trivial $G$-torsor}). If $X = \Spec R$ is affine, we refer to
a $G$-torsor over $X$ also as a $G$-torsor over $R$; it is also
affine because a finite \'etale cover of an affine scheme is affine.
\end{defn}

\begin{defn}
Let $X$ be a $p$-scheme, and let $E \to X$ be a
finite \'etale cover. An \emph{AS-tower for $E \to X$} 
(for ``Artin-Schreier'')
is a sequence of finite \'etale covers
\[
E = E_d \to E_{d-1} \to \cdots \to E_1 \to E_0 = X
\]
in which $E_i \to E_{i-1}$ is equipped with a
$\ZZ/p\ZZ$-torsor structure for $i=1, \dots, d$. 
{}From the transitivity of $p$-typicality (Lemma~\ref{L:ptyptrans}),
we see that the existence of an AS-tower for $E\to X$
implies that $E\to X$ is $p$-typical.
If $X = \Spec R$ and $E = \Spec S$, 
we typically write the tower ring-theoretically,
as $S_0 = R \subset S_1 \subset \cdots \subset S_d = S$,
in which $E_i = \Spec S_i$ and $S_i/S_{i-1}$ is a 
$\ZZ/p\ZZ$-torsor for $i=1,\dots,d$.
\end{defn}

\begin{prop} \label{P:astower}
Let $X$ be a connected $p$-scheme, and
let $E \to X$ be a connected finite \'etale cover.
Then $E\to X$ is $p$-typical if and only if
there exists an AS-tower for $E\to X$.
\end{prop}
\begin{proof}
We have noted already that if there exists an AS-tower for $E \to X$,
then $E \to X$ is $p$-typical (with no connectedness hypotheses).
Conversely, suppose that $E \to X$ is $p$-typical with
monodromy group $G$, which we may assume is nontrivial.
Pick a geometric point $\overline{x}$ of $X$, identify $G$
with the image of $\pi_1(X, \overline{x})$ in
$\Aut(E_{\overline{x}})$, and pick a geometric point
$\overline{y}$ of $E_{\overline{x}}$.
Then the stabilizer of $\overline{y}$ is a proper subgroup of $G$;
thus it is contained in a maximal proper subgroup $H$ of $G$,
which is necessarily normal of index $p$.
In particular, because $H$ is normal, it contains
the stabilizers of \emph{all} of the points of $E_{\overline{x}}$.
Thus $G/H$ is the monodromy
group of a connected $\ZZ/p\ZZ$-torsor $E' \to X$ through which $E$ factors.
By induction, the desired result follows.
\end{proof}

When $E \to X$ is Galois, one gets a bit more.
\begin{prop} \label{P:astower2}
Let $E \to X$ be a connected Galois $p$-typical cover.
Then there exists an AS-tower
$E = E_d \to E_{d-1} \to \cdots \to E_1 \to E_0 = X$
 in which $E_i \to X$ is Galois for $i=1,\dots,d$.
\end{prop}
\begin{proof}
Put $G = \Aut(E \to X)$, which coincides with the monodromy group
of $E \to X$ because the cover is Galois, and assume $G$ is nontrivial.
Since the center of $G$ is nontrivial,
it contains a subgroup $H$ of order $p$, which is normal in $G$.
Let $E_{d-1}$ be the
maximal subcover fixed by $H$, and repeat.
\end{proof}

For $\ZZ/p\ZZ$-torsors over $p$-rings,
one has the following standard result.
\begin{defn}
For $R$ a $p$-ring, a $\ZZ/p\ZZ$-torsor of the form
$S = R[z]/(z^p - z - a)$, in which $1 \in \ZZ/p\ZZ$ acts via $z \mapsto
z + 1$, is called an \emph{Artin-Schreier extension},
or an \emph{AS-extension}, of $R$.
\end{defn}

\begin{prop} \label{P:as}
For any $p$-ring $R$, every $\ZZ/p\ZZ$-torsor of $R$
is an AS-extension. Moreover, two such
torsors $R[z_1]/(z_1^p - z_1 - a_1)$ and $R[z_2]/(z_2^p - z_2 - a_2)$
are isomorphic if and only if $a_1 - a_2 = y^p - y$ for some $y \in R$.
\end{prop}
\begin{proof}
The argument amounts to calculating \'etale cohomology of the
sequence of sheaves:
\[
0 \to \ZZ/p\ZZ \to \GG_a \stackrel{F-1}{\to} \GG_a \to 0
\]
See
\cite{sga4}*{X.3.5},
\cite{katz-canonical}*{1.4.5},
or \cite{milne}*{Proposition~III.4.12} and subsequent discussion.
\end{proof}

\subsection{Connected components in positive characteristic}

We will need to keep careful track of the connected components of certain
AS-towers. Before explaining how we do so, we first make some observations
for arbitrary rings of positive (prime) characteristic.

\begin{lemma} \label{L:p-idem}
Let $R$ be a $p$-ring. Then the set $S = \{x \in R: x^p = x\}$
is the $\FF_p$-subalgebra of $R$ generated by the idempotents
of $R$.
\end{lemma}
\begin{proof}
A straightforward exercise in algebra; alternatively, one may proceed as in
Proposition~\ref{P:as}.
\end{proof}

Counting connected components of rings is closely related to
testing for isomorphisms between finite flat ring extensions, as follows.
\begin{remark} \label{R:isomorphisms}
Let $R$ be a connected $p$-ring,
let $S_1, S_2$ be two connected finite flat
extensions of $R$,
and let $f: S_1 \to S_2$ be an $R$-algebra homomorphism.
Then the graph $\Gamma$ of $f$ is a closed subscheme of $\Spec S_1 \times_R 
\Spec S_2$ which maps isomorphically onto $\Spec S_2$ via the second
projection. In particular, $\Gamma$ is a connected component of
$\Spec S_1 \times_R \Spec S_2$. Conversely, each connected component 
$\Gamma$ of $\Spec S_1 \times_R 
\Spec S_2$ which maps isomorphically onto $\Spec S_2$ via the second projection
corresponds to an $R$-algebra homomorphism $S_1 \to S_2$.
As a consequence, if $g: R \to R'$ is a
ring homomorphism and the induced map
$S_1 \otimes_R S_2 \stackrel{g}{\to}
(S_1 \otimes_R S_2) \otimes_R R'$
induces a bijection of idempotents,
then the induced map
\[
\Hom_{R-\alg}(S_1,S_2) \stackrel{g}{\to} \Hom_{R'-\alg}(S_1 \otimes_R R',S_2
\otimes_R R')
\]
is a bijection.
\end{remark}

\subsection{$p$-injections and $p$-surjections}

We now consider some homomorphisms which behave nicely with respect 
to $p$-typical covers.
\begin{defn}
Let $f: R \to R'$ be a homomorphism of $p$-rings,
and let 
$F$ and $F'$
denote the $p$-power Frobenius maps on $R$ and $R'$,
respectively.
We say $f$ is \emph{$p$-injective} (resp.\ \emph{$p$-surjective})
if the induced functor from $\ZZ/p\ZZ$-torsors over $R$ to
$\ZZ/p\ZZ$-torsors over $R'$ is fully faithful (resp.\ essentially surjective).
These definitions can be reformulated as follows.
\begin{itemize}
\item
The map $f$ is $p$-injective if and only if
$\ker(F-1) \stackrel{f}{\to} \ker(F'-1)$
is surjective and $\coker(F-1) \stackrel{f}{\to} \coker(F'-1)$
is injective.
\item
The map $f$ is $p$-surjective if and only if
$\coker(F-1) \stackrel{f}{\to} \coker(F'-1)$
is surjective.
\end{itemize}
(See the proof of Proposition~\ref{P:p-inj surj} for the explanation of 
how this reformulation follows from Artin-Schreier theory; alternatively,
one may take the reformulation itself as the definition until 
Proposition~\ref{P:p-inj surj} has been proved.)
Using the snake lemma, we may give a second reformulation.
\begin{itemize}
\item
The map $f$ is $p$-injective if and only if
$\ker(f) \stackrel{F-1}{\to} \ker(f)$ is surjective and
$\coker(f) \stackrel{F'-1}{\to} \coker(f)$ is injective
\item
The map $f$ is $p$-surjective if and only if
$\coker(f) \stackrel{F'-1}{\to} \coker(f)$ is surjective.
\end{itemize}
\end{defn}

\begin{remark}
Note that the property of a morphism being $p$-surjective
is not stable under flat base change. For instance,
if $f: R \to R'$ is $p$-surjective but not surjective, then the induced
homomorphism $R[t] \to R'[t]$ is not $p$-surjective.
However, base changing by a $p$-typical extension causes no problems: see
Corollary~\ref{C:p-surjective extends} below.
\end{remark}

\begin{lemma} \label{L:p-surjective extends}
Let $f: R \to R'$ be a homomorphism of $p$-rings,
let
$S = R[z]/(z^p - z - a)$ be an AS-extension of $R$, put
$S' = S \otimes_R R'$, and let $f_S: S \to S'$ be the homomorphism
induced by $f$.
\begin{enumerate}
\item[(a)]
If $f$ induces an injection on idempotents, then so does $f_S$.
\item[(b)]
If $f$ is $p$-injective, then so is $f_S$.
\item[(c)]
If $f$ is $p$-surjective, then so is $f_S$.
\end{enumerate}
\end{lemma}
\begin{proof}
For $l=-1, \dots, p-1$,
let $S_l$ and $S'_l$ be the subsets of $S$ and $S'$, respectively,
consisting of polynomials in $z$ of degree at most $l$ (so that
$S_{-1} = S'_{-1} = \{0\}$); note that each $S_l$ (resp.\ $S'_l$)
is preserved by $F$ (resp.\ by $F'$).
Let $f_l: S_l \to S'_l$ denote the map induced by $f$.
We then have a commutative diagram
\[
\xymatrix{
0 \ar[r] & S_{l-1} \ar[r] \ar^{f_{l-1}}[d] & S_l \ar[r]
\ar^{f_l}[d] & R \ar[r] \ar^{f}[d] & 0 \\
0 \ar[r] & S'_{l-1} \ar[r] & S'_l \ar[r] & R' \ar[r] & 0
}
\]
which by the snake lemma gives rise to a long exact sequence
\[
0 \to \ker(f_{l-1}) \to \ker(f_l) \to \ker(f) \to
\coker(f_{l-1}) \to \coker(f_l) \to \coker(f) \to 0.
\]
We now consider the cases separately.
\begin{enumerate}
\item[(a)]
By Lemma~\ref{L:p-idem} and diagram chasing,
$f$ induces an injection on idempotents if and only if
$F-1$ induces an injection on $\ker(f)$.
In this case, by induction on $l$ and the five lemma,
$F-1$ induces an injection on $\ker(f_l)$ for $l=0,\dots,p-1$.
Taking $l=p-1$, we deduce that $f_S$ induces an injection on
idempotents.
\item[(b)]
If $f$ is $p$-injective, then $F-1$ is surjective on $\ker(f)$
and $F'-1$ is injective on $\coker(f)$. By induction on $l$ and the
five lemma, $F-1$ is surjective on $\ker(f_l)$ and
$F'-1$ is injective on $\coker(f_l)$ for $l=0, \dots, p-1$.
Taking $l=p-1$, we deduce that $f_S$ is $p$-injective.
\item[(c)]
If $f$ is $p$-surjective, then $F'-1$ is surjective on $\coker(f)$.
By induction on $l$ and the five
lemma, $F'-1$ is surjective on $\coker(f_l)$ for $l=0, \dots, p-1$.
Taking $l=p-1$, we deduce that $f_S$ is $p$-surjective.
\end{enumerate}
\end{proof}
\begin{cor} \label{C:p-surjective extends}
Let $f: R \to R'$ be a homomorphism of $p$-rings,
let $R = S_0 \subset S_1 \subset \cdots \subset S_d$ be
an AS-tower over $R$, put $S'_i = S_i \otimes_R R'$ for $i=1,\dots, d$,
and let $f_i: S_i \to S'_i$ be the homomorphism induced by $f$.
\begin{enumerate}
\item[(a)]
If $f$ induces an injection on idempotents, then so does each $f_i$.
\item[(b)]
If $f$ is $p$-injective, then so is each $f_i$.
\item[(c)]
If $f$ is $p$-surjective, then so is each $f_i$.
\end{enumerate}
\end{cor}

\begin{prop} \label{P:p-inj surj}
Let $f: R \to R'$ be a homomorphism of $p$-rings.
Let $\calS_R$ and $\calS_{R'}$ be the categories of AS-towers over $R$
and $R'$, respectively, in which the only morphisms are isomorphisms of towers.
\begin{enumerate}
\item[(a)]
If the map 
$f$ is $p$-injective, then the base change functor $f^*: \calS_{R} \to 
\calS_{R'}$
is fully faithful.
\item[(b)]
The map $f$ is $p$-surjective if and only if
the base change functor $f^*: \calS_R \to \calS_{R'}$ 
is essentially surjective.
\end{enumerate}
\end{prop}
\begin{proof}
\begin{enumerate}
\item[(a)]
Suppose that $f$ is $p$-injective.
Given two AS-towers
$R = S_0 \subset S_1 \subset \cdots \subset S_d = S$
and
$R = T_0 \subset T_1 \subset \cdots \subset T_d = S$
which become isomorphic over $R'$,
write $S_1 = R[y]/(y^p - y - a)$ and $T_1 = R[z]/(z^p - z - b)$.
By Proposition~\ref{P:as}, $f(a)$ and $f(b)$ represent the same element
of $\coker(F'-1)$; hence they also represent the same element of
$\coker(F-1)$. Thus $S_1 \cong T_1$; moreover,
by Lemma~\ref{L:p-surjective extends}, the map $S_1 \to S_1 \otimes_R R'$
is $p$-injective. Repeating the argument, we see that the two towers
are isomorphic, and so $f^*$ is fully faithful.

\item[(b)]
Suppose that $f^*$ is essentially surjective. 
Let $S' = R'[z]/(z^p - z - a)$ be an AS-extension of $R'$; by hypothesis,
there exists an AS-extension $S = R[z]/(z^p - z - b)$ such that
$S \otimes_R R' \cong S'$ as a $\ZZ/p\ZZ$-torsor. By Proposition~\ref{P:as},
we must have $f(b) - a = y^p - y$ for some $y \in R'$. We deduce
that the map $\coker(F-1) \to \coker(F'-1)$ induced by $f$ is surjective,
and so $f$ is $p$-surjective.

Conversely, suppose that $f$ is $p$-surjective.
Given an AS-tower
$R' = S'_0 \subset S'_1 \subset \cdots \subset S'_d = S'$,
we construct a corresponding AS-tower
$R = S_0 \subset S_1 \subset \cdots \subset S_d$ inductively as follows.
Start with $S_0 = R$. Given $S_0, \dots, S_i$
and an isomorphism $S_i \otimes_R R' \cong S'_i$,
note that $f: S_i \to S'_i$ is $p$-surjective by
Lemma~\ref{L:p-surjective extends}. By
Proposition~\ref{P:as}, we can write $S'_{i+1} = S'_i[z]/(z^p - z - a)$
for some $a \in f(S_i)$; we may then set $S_{i+1} = S_i[z]/(z^p - z - b)$
for any $b \in S_i$ with $f(b) = a$. Thus the inductive construction
continues, and so $f^*$ is essentially surjective.
\end{enumerate}
\end{proof}

\begin{remark}
Beware that proving results about the category of AS-towers does not
immediately yield results about $p$-typical covers; for that, stronger
connected hypotheses are needed, as in the next section.
\end{remark}

\subsection{$p$-limits and canonical extensions}

\begin{convention}
Given a partially ordered set $S$, we view $S$ as a category in
which $\Mor(s,t)$ is a singleton set if $s \geq t$ and is empty
otherwise.
\end{convention}

\begin{defn}
A \emph{diagram} in a category $\calC$
is a functor $D$ from a partially ordered set $S$ 
to $\calC$; we call $S$ the \emph{support} of $D$.
Given a subset $T$ of $S$, let $D_{T}$ denote the restriction
of $D$ to $T$.
\end{defn}

\begin{defn}
Given a diagram $D$ with support $S$,
put $S_1 = S_2 = S \cup \{s'\}$ for some $s' \notin S$, 
and extend the partial order from $S$ to 
$S_1$ and $S_2$ by declaring that in $S_1$,
$s' \geq s$ and $s \not\geq s'$ for all $s \in S$,
while in $S_2$, $s \geq s'$ and $s' \not\geq s$ for all $s \in S$.
For an object $X \in \calC$, a \emph{morphism from $X$
to $D$} (resp.\ a \emph{morphism from $D$ to $X$})
is a diagram $D'$ supported on $S_1$ (resp.\ on $S_2$) with 
$D'(s') = X$ and $D'_S = D$; let $\Mor(X,D)$ (resp.\ $\Mor(D,X)$)
denote the set of these morphisms.
A \emph{limit} (resp.\ \emph{colimit})
of a diagram $D$ is an object $X \in \calC$ representing 
the functor $Y \mapsto \Mor(Y,D)$ (resp.\ the functor
$Y \mapsto \Mor(D,Y)$); by construction, a (co)limit is unique
up to unique isomorphism if it exists.
\end{defn}

\begin{remark}
Note that every diagram in the category of affine $p$-schemes
has a limit, which can be constructed by repeatedly constructing
products and equalizers. (Arbitrary products are given by ``infinite
tensor products'', which are generated by terms
which have the factor 1 in all but finitely many places.)
However, a diagram in the category of \emph{connected}
affine $p$-schemes need not have a limit.
\end{remark}

\begin{defn}
Let $D$ be a nonempty diagram in the category of connected
affine $p$-schemes.
A \emph{$p$-limit} of $D$ is a
connected affine $p$-scheme $Y$ equipped with a morphism
$Y \to D$, which becomes a colimit of $D$ in the category of abelian groups
upon applying the contravariant functor $X \mapsto 
\coker(F-1, \Gamma(X,\calO_X))$.
\end{defn}
\begin{remark} \label{R:p-inj surj}
Note that $D$ admits a limit $X$ in the category of affine $p$-schemes,
and that if $Y$ is a $p$-limit of $D$, then the induced 
homomorphism $\Gamma(X,\calO_X) \to \Gamma(Y,\calO_Y)$ is $p$-surjective:
the direct sum of the $\coker(F-1,\Gamma(Z,\calO_Z))$ 
for all $Z$ in the diagram surjects onto $\coker(F-1,\Gamma(Y,\calO_Y))$,
but this surjection factors through $\coker(F-1,\Gamma(X,\calO_X))$.
\end{remark}

\begin{theorem} \label{T:p-limit}
Let $D$ be a nonempty finite diagram in the category of connected
affine $p$-schemes, let $S$ be the support of $D$,
and let $Y \to D$ be a morphism; for $s \in S$, let $f_s$ be the induced
morphism
from $Y$ to $D(s)$.
Choose a geometric point $\overline{y}$ of $Y$. Then 
$Y$ is a $p$-limit of $D$ if and only if
$\pi_1^p(Y, \overline{y})$ is a limit, in the category of pro-$p$-groups,
of the diagram induced by $D$ on the $\pi_1^p(D(s), f_s(\overline{y}))$.
\end{theorem}
\begin{proof}
First suppose that $Y$ is a $p$-limit of $D$. Then the homomorphism
\begin{equation} \label{eq:pi1p}
\pi_1^p(Y, \overline{y}) \to \lim \pi_1^p(D(s), f_s(\overline{y}))
\end{equation}
is seen to be injective as follows. Given a non-identity element $\tau$
of $\pi_1^p(Y, \overline{y})$, choose a connected $p$-typical cover
$E$ of $Y$ such that $\tau$ acts nontrivially on $E_{\overline{y}}$.
By Proposition~\ref{P:astower}, $E$ admits an AS-tower;
by Remark~\ref{R:p-inj surj} and Proposition~\ref{P:p-inj surj}, that AS-tower
can be obtained by pullback from some AS-tower over a limit of
$D$ in the category of affine $p$-schemes. Hence 
the image of $\tau$ in $\lim \pi_1^p(D(s), f_s(\overline{y}))$
is not the identity element, so \eqref{eq:pi1p} is injective.

Suppose now that \eqref{eq:pi1p} fails to be surjective.
We can then construct a nontrivial
continuous homomorphism $g:
\lim \pi_1^p(D(s), f_s(\overline{y})) \to \ZZ/p\ZZ$
(for the discrete topology on $\ZZ/p\ZZ$) whose restriction
to $\pi_1^p(Y, \overline{y})$ is trivial. 
For $s \in S$, put $C_s = \coker(F-1, \Gamma(D(S),\calO))$.
We then obtain from $g$ and Proposition~\ref{P:as}
an element $c_s \in C_s$ for each $s \in S$,
such that if $s \to t$ is a morphism in $S$, then the corresponding
morphism $C_t \to C_s$ carries $c_t$ to $c_s$.
Since $Y$ is a $p$-limit, the $c_s$ correspond to a nonzero element
of $\coker(F-1, \Gamma(Y,\calO_Y))$,
which gives rise to a nontrivial $\ZZ/p\ZZ$-torsor on $Y$. This contradicts
the fact that $g$ restricts trivially to $\pi_1^p(Y, \overline{y})$;
the contradiction yields the surjectivity of \eqref{eq:pi1p}, as desired.

We have now shown that if $Y$ is a $p$-limit of $D$, then
\eqref{eq:pi1p} is an isomorphism. Suppose now conversely that
\eqref{eq:pi1p} is an isomorphism. Then
the maximal elementary abelian quotient of
$\pi_1^p(Y, \overline{y})$ is the limit, in the category of elementary
abelian $p$-groups, of the maximal elementary abelian quotients of the
$\pi_1^p(D(s), f_s(\overline{y}))$. But by Proposition~\ref{P:as},
these quotients are dual to the cokernels of $F-1$ on these schemes.
Hence $Y$ is a $p$-limit of $D$, as desired.
\end{proof}

Theorem~\ref{T:p-limit} may be a bit obscure as written; some of its
corollaries may be more edifying.
\begin{defn}
Let $f: R \to R'$ be a homomorphism of connected $p$-rings,
and let $F$ and $F'$ be the $p$-power Frobenius maps on $R$ and $R'$,
respectively. We say $f$ is \emph{$p$-faithful} if
the induced map $\coker(F-1) \stackrel{f}{\to} \coker(F'-1)$ 
is a bijection.
\end{defn}
\begin{cor} \label{C:p-faithful}
Let $f: R \to R'$ be a homomorphism of connected $p$-rings,
choose a geometric point $\overline{x}'$ of $\Spec R'$,
and put $\overline{x} = f(\overline{x}')$.
Then $f$ is $p$-faithful if and only if $\pi_1^p(\Spec R', \overline{x}')
\stackrel{f}{\to} \pi_1^p(\Spec R, \overline{x})$ is a bijection.
\end{cor}
\begin{example} \label{E:katz example}
For any $p$-ring $R$, the canonical inclusion $f: R[t^{-1}] \to R((t))$
is $p$-faithful: the kernel of $f$ is trivial, and the cokernel of $f$
is isomorphic as a Frobenius module to
$t R \llbracket t \rrbracket$, on which $F-1$ is bijective.
The conclusion of Corollary~\ref{C:p-faithful} in this case is a result
of Katz \cite{katz-canonical}*{Proposition~1.4.2}.
Although Katz's proof looks different (it involves manipulating
the cohomology of pro-$p$-groups), our proof is basically a transcription
of Katz's argument into the language of AS-towers.
\end{example}

\begin{example} \label{E:Laurent series}
Let $R$ be a $p$-ring.
Consider the diagram consisting of the two natural maps
$\Spec R[t] \to \Spec R$ and $\Spec R[t^{-1}] \to \Spec R$.
Then $\Spec R[t,t^{-1}]$ is a $p$-limit of this diagram; we thus have an
isomorphism
\[
\pi_1^p (\Spec R[t,t^{-1}]) \to
\pi_1^p(\Spec R[t]) \times_{\pi_1^p(\Spec R)} \pi_1^p(\Spec R[t^{-1}])
\]
after choosing basepoints. (Namely, choose a geometric point of
$\Spec R[t,t^{-1}]$ and obtain the other basepoints by applying the
maps in the diagram.)
\end{example}

Here is a slight variation of the previous example.
\begin{cor}
Let $R$ be an $\overline{\FF_p}$-algebra. Then every $p$-typical extension
of $R[t]$ is contained in the tensor product of a $p$-typical
extension of $R[t]$ in which $R$ is integrally closed, and a
$p$-typical extension of $R[t]$ obtained by base change from $R$.
\end{cor}
\begin{proof}
Put $R' = \overline{\FF_p} + t R[t] \subseteq R[t]$. Then
$\Spec R[t]$ is a $p$-limit of the diagram consisting of
$\Spec R$ and $\Spec R'$ with no arrows, so by
Theorem~\ref{T:p-limit}, we have
$\pi_1^p(\Spec R[t]) \cong \pi_1^p(\Spec R) \times 
\pi_1^p(\Spec R')$
(for appropriate basepoints). Thus every $p$-typical extension of $R[t]$
is contained in the tensor product of a $p$-typical extension obtained
by base change from $R$,
and a $p$-typical extension obtained by base change
from $\overline{\FF_p} + t R[t]$;
in the latter, the restriction to the $t=0$ locus must split completely,
so $R$ must be integrally closed. This yields the desired result.
\end{proof}
\begin{remark}
This corollary should be a bit surprising: for a general finite \'etale
extension of $R[t]$, or even of $R((t^{-1}))$, one cannot split off the
residual extension in this fashion. For instance, if the extension
is obtained by adjoining $z$ with $z^p - z = at$ for $a$ in some
finite \'etale extension of $R$, it is typically impossible to present
the extension as in the corollary \emph{unless} $a$ generates a
$p$-typical extension of $R$ (in which case the corollary applies).
\end{remark}

\section{Complexity measures for $p$-typical extensions}

We next propose a mechanism for handling the ``complexity''
of a $p$-typical extension, via what we call ``height functions'';
the mechanism is modeled on basic ramification theory for complete
discretely valued fields. As in other instances where complexity-bounding
functions arise (e.g., Diophantine approximation, from which the term
``height function'' was borrowed), it is a bit tedious to introduce and
deal with such functions, but things are made a bit easier by the fact that
the intended use of these functions permits one to be somewhat sloppy in
dealing with them. The reader impatient to get to some meaningful results
may wish to skip ahead to the next chapter before continuing here.

\subsection{Ramification filtrations for local fields}

The model for our height functions
is the highest break function coming from the
ramification filtration on the Galois group of $k((t))$, so we start
by reviewing that construction.
For all unproved assertions in this section, see \cite{serre}.

\begin{defn} \label{D:lower}
Let $F$ be a complete discretely valued field whose residue field $k$
is perfect (e.g., the power series field $k((t))$).
Let $E/F$ be a finite Galois field extension with group $G$,
let $\gotho_E$ and $\gotho_F$ be the valuation subrings of $E$ and $F$,
and let $v_E$ be the valuation on $E$,
normalized so that $v_E$ maps $E^*$ onto $\ZZ$.
For $i \geq -1$, let $G_i$ be the subgroup of $g \in G$
for which $v_E(a^g - a) \geq i+1$ for all $a \in \gotho_E$;
the decreasing filtration $\{G_i\}$ is called 
the \emph{lower numbering filtration} of $G$ \cite{serre}*{\S IV.1}.
\end{defn}

\begin{defn} \label{D:upper}
With notation as in Definition~\ref{D:lower}, 
define the function
\[
\phi_{E/F}(u) = \int_{0}^u \frac{dt}{[G_0:G_t]}.
\]
Then $\phi_{E/F}$ is a homeomorphism of $[-1, \infty)$ with itself;
let $\psi_{E/F}$ denote the inverse function.
Define the \emph{upper
numbering filtration} of $G$ by $G^i = G_{\psi_{E/F}(i)}$
\cite{serre}*{\S IV.3}. It has the property that 
if $E'/F$ is a Galois subextension of $E/F$ with Galois group $H$, then
the image of each $G^i$ under the natural surjection $G \twoheadrightarrow H$
is precisely $H^i$; this follows from Herbrand's theorem
\cite{serre}*{Proposition~IV.14}.
\end{defn}

\begin{defn} \label{D:highest break}
For $F$ as in Definition~\ref{D:lower}
and $E/F$ a finite Galois field extension,
define the \emph{highest break} of $E$,
denoted $b(E/F)$,
to be the largest $i$ such that $G^i \neq G^j$ for any $j>i$, or zero
if no such $i$ exists. 
If $E/F$ is a field extension which is finite separable but not Galois,
we define $b(E/F) = b(E'/F)$, for $E'/F$ the Galois closure of $E/F$.
If $E$ is not a field but only an \'etale $F$-algebra, we define
$b(E/F)$ to be the maximum highest break of any component of $E$.
With these rules, one has the following properties.
\begin{enumerate}
\item[(a)]
$b(F/F) = 0$ (evident).
\item[(b)]
If $E'$ is an $F$-subalgebra of $E$, then $b(E'/F) \leq b(E/F)$ (evident).
\item[(c)]
$b((E_1 \oplus E_2)/F) = \max\{b(E_1/F), b(E_2/F)\}$ (formal).
\item[(d)]
$b((E_1 \otimes E_2)/F) = \max\{b(E_1/F), b(E_2/F)\}$ (not formal,
but follows from Herbrand's theorem).
\item[(e)]
If $E/F$ is a Galois field extension and $E'$ is an \'etale $E$-algebra, then
$b(E'/F) = \max\{b(E/F), \phi_{E/F}(b(E'/E))\}$ (because the lower numbering
is stable under taking subgroups).
\end{enumerate}
\end{defn}

\begin{remark} \label{R:imperfect}
In case $k$ is imperfect, there are several competing analogues of the
upper numbering filtration; these include the ``residual perfection''
construction of Borger \cite{borger}, and the ``nonlogarithmic''
and ``logarithmic'' rigid geometric constructions of Abbes and Saito
\cite{abbes-saito}. We will not use any of these in this paper.
\end{remark}

\subsection{Artin-Schreier extensions and highest breaks}

We next recall some standard facts about Artin-Schreier extensions of
a power series field.
\begin{lemma}  \label{L:as1}
For $k$ a perfect $p$-field, and
for $a \in F = k((t))$, put
\[
m = \inf_{x \in F} \{-v_{F}(a - x^p + x) \}
\]
and put $E = F[z]/(z^p - z - a)$. Then the following hold.
\begin{enumerate}
\item[(a)] Either $m = -\infty$ (that is, $E$ is not a field) or $m \geq 0$.
\item[(b)] If $m \geq 0$, then the extension $E/F$ is unramified if and only if
$m = 0$.
\item[(c)] If $m > 0$, then $m$ is not divisible by $p$,
and $E/F$ has highest break $m$.
\end{enumerate}
\end{lemma}
\begin{proof}
\begin{enumerate}
\item[(a)]
Suppose $m<0$, which means that there exists $y \in F$ such that
$b = y^p - y - a$ satisfies $v_F(b) > 0$. Then the series
$b + b^p + b^{p^2} + \cdots$ converges in $F$, and its limit $c$ satisfies
$c - c^p = b$. This yields $a = (c+y)^p - (c+y)$, and so $m=-\infty$.
\item[(b)]
Note that $m=0$ implies, as in (a), that $y^p - y - a \in k$
for some $y \in F$, and so $E$ is unramified. Conversely, if $E$ is unramified,
then the residue field
$\overline{E}$ of $E$ must be an Artin-Schreier extension of $k$, say
$k[y]/(y^p - y - b)$ with $b \in k$. If we choose $b$ so that the
$\ZZ/p\ZZ$-torsor structures on $E/F$ and $\overline{E}/k$
are compatible, by Proposition~\ref{P:as} we must then have $a - b = x^p - x$
for some $x \in F$, yielding $m = 0$.
\item[(c)]
If $a = c_{-pn} t^{-pn} + \cdots$, then $a - c_{-pn} t^{-pn} + c_{-pn}^{1/p}
t^{-n}$ has strictly larger valuation than does $a$. 
Hence if $m$ is positive, it cannot
be divisible by $p$. To compute the highest break, pick integers $r,s$
with $r>0$ and $-rm + sp = 1$, and put $u = z^r t^s$; then $v_E(u) = 1$, i.e.,
$u$ is a uniformizer of $E$. By \cite{serre}*{Proposition IV.5},
the highest break of $E/F$ equals $v_E(u'/u - 1)$, where $u'$ is the
image of $u$ under the automorphism $z \mapsto z+1$ of $E$.
Since $r$ is not divisible by $p$, we have
\[
\frac{u'}{u} = (z+1)^r z^{-r} = 1 + rz^{-1} + \cdots
\]
and so $v_E(u'/u - 1) = v_E(z^{-1}) = m$, as desired.
\end{enumerate}
\end{proof}

One can also obtain a bound on the highest break in an AS-tower.
The following bound is not optimal, but it suffices for our purposes.
\begin{cor} \label{C:asbound}
Let $k$ be a perfect $p$-field,
let $k((t)) = E_0 \subset E_1 \subset \cdots \subset E_d = E$ be
an AS-tower, and choose $\ell \geq 1$ such that 
for $i=1, \dots, m$, $E_i \cong E_{i-1}[z]/(z^p - z - c_i)$ for some $c_i$ with
$v_{E_0}(c_i) \geq -\ell$. Then $b(E/E_0) \leq d\ell$.
\end{cor}
\begin{proof}
We proceed by induction on $d$, the case $d=1$ following from
Lemma~\ref{L:as1}. For $d>1$, if $E_1/E_0$ is disconnected, then
we can correspondingly split $E$ as a direct sum $E'_1 \oplus \cdots
\oplus E'_p$, in which for $j=1, \dots, p$, $E'_j$ admits an AS-tower
of length $d-1$ over $E_0$. By the induction hypothesis, we have
$b(E/E_0) = \max_j \{b(E'_j/E_0)\} \leq (d-1) \ell$.

If $E_1/E_0$ is connected, by Lemma~\ref{L:as1}
we have $b(E_1/E_0) = m$ for some nonnegative integer $m \leq \ell$,
and by the induction hypothesis we have
 $b(E/E_1) \leq (d-1)(p\ell)$.
For $x \geq m$,
$\phi_{E_1/E_0}(x) = m + (x-m)/p$, so
\begin{align*}
b(E/E_0) &\leq \phi_{E_1/E_0}((d-1)p\ell ) \\
&= m + \frac{(d-1)p\ell - m}{p} \\
&= (d-1)\ell + \frac{m(p-1)}{p} \\
&\leq d\ell,
\end{align*}
as desired.
\end{proof}

We also need to know that the highest break drops under
specialization.
\begin{prop} \label{P:Laumon}
Let $R \to R'$ be a surjective morphism of perfect $p$-domains,
let $S$ be a $p$-typical extension of $R((t))$, and put
$S' = S \otimes_{R((t))} R'((t))$. Let $K$ and $K'$ be the fraction
fields of $R$ and $R'$, respectively. Then
\[
b(S \otimes_{R((t))} K((t))/K((t))) \geq 
b(S' \otimes_{R'((t))} K'((t))/K'((t))).
\]
\end{prop}
\begin{proof}
This follows from the Deligne-Laumon semicontinuity theorem
\cite{laumon}.
\end{proof}

\subsection{Presentations of AS-towers}

To talk about height functions on $p$-typical extensions of more general
rings, we need to fix a bit of terminology concerning presentations
of AS-towers.

\begin{defn}
Given an AS-tower $R = S_0 \subset S_1 \subset \cdots \subset S_d = S$
over a $p$-ring $R$,
a \emph{presentation of $S$} is a sequence of isomorphisms
\[
S_i \cong S_{i-1}[z_i]/(z_i^p - z_i - 
P_i(z_1,\dots,z_{i-1})) \qquad (i=1,\dots,d),
\]
where $P_i(z_1, \dots, z_{i-1})$ is a polynomial over $R$
of degree at most $p-1$ in each variable;
by Proposition~\ref{P:as},
such a presentation always exists.
Given a presentation of $S$, each element $x \in S$ can be written
uniquely as a polynomial in $z_1,\dots,z_d$ over $R$ with
degree at most $p-1$ in each variable; we call this polynomial the
\emph{minimal representation} of $x$.
\end{defn}

In terms of presentations, one has the following evident but useful lemma.
\begin{lemma} \label{L:degrees}
Given an AS-tower $R = S_0 \subset S_1 \subset \cdots \subset S_d = S$
over a $p$-ring $R$, and a presentation
\[
S_i \cong S_{i-1}[z_i]/(z_i^p - z_i - 
P_i(z_1,\dots,z_{i-1})) \qquad (i=1,\dots,d)
\]
of $S$, choose integers $j_1, \dots, j_d \in \{0, \dots, p-1\}$,
and put $x = z_1^{j_1} \cdots z_d^{j_d}$. Then the minimal representation
of $x^p$, written as a polynomial in $z_d$ over $S_{d-1}$, 
is monic of degree $j_d$.
\end{lemma}
\begin{proof}
Note that for each $i$, $z_i^p$ can be rewritten as $z_i$ plus a polynomial
in the preceding variables; this implies the claim.
\end{proof}

\begin{defn}
If $V$ is an additive subgroup of $R$, we say a presentation of $S$
is \emph{defined over $V$} if each $P_i$ has its coefficients in $V$.
\end{defn}

\subsection{Height functions}

\begin{defn} \label{D:height function}
Let $R_0$ be a connected $p$-ring, and let
$R$ be a connected $R_0$-algebra.
A \emph{height function (over $R_0$)} on $\calC^p_R$ (the category of
$p$-typical extensions of $R$) is a function $h$ from the set of
isomorphism classes of elements of $\calC^p_R$ to the nonnegative real numbers,
having the following properties.
\begin{enumerate}
\item[(a)] $h(S_1 \oplus S_2)$ is bounded above by some function of
$h(S_1), h(S_2), \deg(S_1/R), \deg(S_2/R)$.
\item[(b)] $h(S_1 \otimes S_2)$ is bounded above by some function
of $h(S_1), h(S_2), \deg(S_1/R), \deg(S_2/R)$.
\item[(c)] If $S_1 \subseteq S_2$, then $h(S_1)$
is bounded above by some function of $h(S_2), \deg(S_2/R)$.
\item[(d)] For any positive integer $d$ and any finite $R_0$-submodule
$V$ of $R$, there exists a nonnegative real number $\ell$ such that
for any connected
AS-tower $R = S_0 \subset S_1 \subset \cdots \subset S_d = S$
admitting a presentation defined over $V$, we have $h(S) \leq \ell$.
\item[(e)] 
For any positive integer $d$ and any nonnegative real number $\ell$,
there exists a finite $R_0$-submodule 
$V$ of $R$ such that for any connected AS-tower 
$R = S_0 \subset S_1 \subset \cdots \subset S_d = S$
with $h(S) \leq \ell$, there exists a presentation of $S$ defined over $V$.
\end{enumerate}
We say
$h$ is a \emph{strong height function} if the following additional
conditions hold.
\begin{enumerate}
\item[(a$'$)] $h(S_1 \oplus S_2) \leq \max\{h(S_1), h(S_2)\}$.
\item[(b$'$)] $h(S_1 \otimes S_2) \leq \max\{h(S_1), h(S_2)\}$.
\item[(c$'$)] If $S_1 \subseteq S_2$, then $h(S_1) \leq h(S_2)$.
\end{enumerate}
We extend a height function to continuous homomorphisms $\rho: \pi^p_1(\Spec R,
\overline{x}) \to G$,
for $\overline{x}$ a geometric point of $\Spec R$ and 
$G$ a finite discrete group, by declaring that $h(\rho) = h(S)$,
where $S \in \calC^p_R$ is chosen so that
$\pi_1^p(\Spec S, \overline{y})$ is the kernel of $\rho$
(for an appropriate geometric point $\overline{y}$ of $\Spec S$).
\end{defn}
\begin{lemma}
With notation as in Definition~\ref{D:height function},
let $R'$ be a connected $p$-typical extension of $R$.
Then any height function $h$ over $R_0$ on $\calC^p_R$ induces a height
function $h'$ over $R_0$ on $\calC^p_{R'}$ (given by $h'(S) = h(S)$).
\end{lemma}
\begin{proof}
Straightforward.
\end{proof}

The condition (e) is not so easy to check directly, but fortunately
one need only verify it for Artin-Schreier extensions, as confirmed
by the following proposition.

\begin{prop} \label{P:height1}
Given conditions (a)-(d) of Definition~\ref{D:height function}, 
if condition (e)
holds for $d=1$, then it holds for all $d$.
\end{prop}
\begin{proof}
We proceed by induction on $d$ (simultaneously for all $\ell$), 
the case $d=1$ being the input hypothesis.
Given the claim for $d-1$ and given a connected AS-tower 
$R=S_0 \subset S_1 \subset \cdots \subset S_d = S$ with $h(S) \leq \ell$,
we may choose
a presentation for $S_{d-1}$ over some
finite $R_0$-module
depending only on $d$ and $\ell$.

Write $S_d = S_{d-1}[z_{d}]/(z_{d}^p - z_{d} - a_{d-1})$
and write $a_{d-1} = \sum_{i=0}^{p-1} c_{i} z_{d-1}^i$
for $c_{i} \in S_{d-2}$.
Let $j$ be the degree of $a_{d-1}$ as a polynomial in $z_{d-1}$,
so that $c_j \neq 0$
but $c_{j+1} = \cdots = c_{p-1} = 0$. We prove 
that for some $w \in S_{d-1}$ of degree at most $j$
as a polynomial in $z_{d-1}$, the coefficients in the minimal representation
of $a_{d-1} - w^p  + w$
lie in some finite $R_0$-module depending only on $d,\ell,j$.
The proof of this claim constitutes an inner induction on $j$.

For $j=0$, we may apply the outer induction hypothesis to 
$S_{d-2}[z]/(z^p - z - c_0)$ to deduce the claim.
Otherwise,
let $g$ be the automorphism of $S_{d-1}$
over $S_{d-2}$ given by $z_{d-1} \mapsto z_{d-1}+1$, and
define the map $\Delta: S_{d-1} \to S_{d-1}$ by
$\Delta(x) = x^g - x$. Then the $j$-th tensor power of $S_d$ over $S_{d-2}$,
which has height bounded by a function of $d,\ell,j$ by (b),
contains 
\[
S_{d-1}[z]/(z^p - z - \Delta^{(j)}(a_{d-1}));
\]
since $\Delta^{(j)} (\sum_{i=0}^j c_i z_{d-1}^i) = j! c_j$, 
we deduce that
$S_{d-2}[z]/(z^p - z - j! c_j)$ has height bounded by a function
of $d,\ell,j$.

Applying the outer induction hypothesis yields $w' \in S_{d-2}$ such that
the minimal representation of $b' = c_j - (w')^p + w'$ has coefficients in some
finite $R_0$-module depending only on $d,\ell,j$. Put
\[
a_{d-1}' = a_{d-1} - (w' z_{d-1}^j)^p + w' z_{d-1}^j,
\]
so that $S_d = S_{d-1}[z]/(z^p - z - a_{d-1}')$.
Then $a_{d-1}' - b' z_{d-1}^j$ has degree at most $j-1$ as a polynomial in
$z_{d-1}$. If we put
\begin{align*}
S' &= S_{d-1}[x]/(x^p - x - b' z_{d-1}^j), \\
S'' &= S_{d-1}[y]/(y^p - y - a_{d-1}' - b' z_{d-1}^j),
\end{align*}
then the height of $S'$ is bounded by some function of $d,\ell,j$
by condition (d) of Definition~\ref{D:height function}.
On the other hand, $S''$ is contained in $S' \otimes_{S_{d-1}} S_d$,
and the heights of $S_d$ and $S'$ are bounded by some function
of $d,\ell,j$, so the same is true of $S''$.
Applying the inner
induction hypothesis to $S''$, we obtain $w'' \in S_{d-1}$ of degree
at most $j-1$ as a polynomial in $z_{d-1}$, such that
$a_{d-1}' - b' z_{d-1}^j - (w'')^p + w''$ has coefficients in some
finite $R_0$-module depending only on $d,\ell,j$.
We may then take $w = w' z_{d-1}^j + w''$, as
\[
a_{d-1} - w^p + w = b' z_{d-1}^j + 
(a_{d-1}' - b' z_{d-1}^j - (w'')^p + w'')
\]
has coefficients in some finite $R_0$-module depending on $d,\ell,j$.
This completes the
proof of the inner induction.

The inner induction for $j=p-1$ implies the outer induction,
so the proof is complete.
\end{proof}

\begin{example} \label{E:onedim}
For $R = k((t))$ with $k$ a perfect $p$-field, 
the highest break function $h(S) = 
b(S/k((t)))$ is a strong height
function on $\calC^p_{k((t))}$ over $k$: 
properties (a$'$), (b$'$), (c$'$) follow from
Definition~\ref{D:highest break},
property (d) from Corollary~\ref{C:asbound}, and property
(e) for $d=1$ from
Lemma~\ref{L:as1}.
\end{example}

\begin{remark} \label{R:onedim}
Already for $R=k((t))$ with $k$ an imperfect $p$-field, 
it is less than evident how to
construct a height function on $\calC^p_{k((t))}$ over $k$, 
since the na\"\i ve highest break function
$b(S \otimes k^{\perf}((t))/k^{\perf}((t)))$ will not do. To see this,
choose $c \in k \setminus k^p$, then note that
the heights of $k((t))[z]/(z^p - z - c t^{-p^n})$ would all be equal to 1,
whereas these extensions do not simultaneously admit presentations defined over
some finite dimensional $k$-vector space.
It should be possible to extract a height function from any of the
constructions of a ramification filtration 
mentioned in Remark~\ref{R:imperfect}, but we have not attempted to do so.
\end{remark}

\section{$p$-typical covers of affine toric varieties}
\label{sec:ramif2}

In this chapter, we study the $p$-typical fundamental groups of
affine toric varieties; while the case of ordinary affine space is doubtless
the most important, it is not any harder to work in this generality.
Our main results in this direction are some decomposition theorems
for these fundamental groups (Theorem~\ref{T:discrete} and its consequence
Theorem~\ref{T:splits2}).

\setcounter{equation}{0}
\begin{convention} \label{conv:alg closed}
Throughout this chapter, let $R$ denote a connected $p$-ring.
\end{convention}

\subsection{Some toric rings}

\begin{defn}
Define a \emph{convex cone} in $\RR^n$ as a nonempty
subset $\sigma \subseteq \RR^n$ such that:
\begin{enumerate}
\item[(a)]
if $\bv \in \sigma$, then $c \bv \in \sigma$ for any $c \in \RR_{\geq 0}$;
\item[(b)]
if $\bv, \bw \in \sigma$, then
$c\bv + (1-c)\bw \in \sigma$ for any $c \in [0,1]$.
\end{enumerate}
Note that the intersection of convex cones is again a convex cone;
we say the convex cone $\sigma$ is \emph{finitely generated} if it can
be written as a finite intersection of open and closed halfspaces.
\end{defn}


\begin{defn}
Given a convex cone $\sigma$, let $R_\sigma$ denote the monoid algebra
$R[\sigma \cap \ZZ^n]$; 
for convex cones $\sigma, \tau$ with
$\sigma \subseteq \tau$, there is a natural inclusion
$R_\sigma \subseteq R_\tau$.
Given an element $x \in R_\sigma$, write
\[
x = \sum_{\bv \in \sigma \cap \ZZ^n} c_{\bv} [\bv],
\]
and define the \emph{support} of $x$ to be the set of 
$\bv \in \sigma \cap \ZZ^n$ such that $c_{\bv} \neq 0$.
\end{defn}
\begin{remark}
If $\sigma$ is a convex cone equal to the intersection of finitely many
closed halfspaces defined by
linear functionals over $\QQ$, then $\Spec R_\sigma$ is an affine
toric variety, and conversely. (Note that in our terminology, toric varieties
are necessarily normal.)
\end{remark}

\begin{convention}
For the rest of the chapter, fix a geometric point
$\overline{x}$ of $\Spec R_{\RR^n}$; we may also view
$\overline{x}$ as a geometric point of $\Spec R_\sigma$ for any
convex cone $\sigma \subseteq \RR^n$. We will thus drop this basepoint
from the notation when considering the fundamental group of $\Spec R_\sigma$.
\end{convention}

\begin{prop} \label{P:split1}
Suppose that 
$\sigma, \sigma_1, \dots, \sigma_n$ are convex cones
with $\sigma = \sigma_1 \cup \cdots \cup \sigma_n$. Then
$\Spec R_\sigma$ is a $p$-limit of the diagram consisting of the arrows
$\Spec R_{\sigma_i} \to \Spec R_{\sigma_i \cap \sigma_j}$
for $1 \leq i,j \leq n$. Consequently,
the group $\pi_1^p(\Spec R_\sigma)$ is a limit
of the diagram consisting of the arrows
$\pi_1^p(\Spec R_{\sigma_i}) \to \pi_1^p(\Spec R_{\sigma_i 
\cap \sigma_j})$
for $1 \leq i,j \leq n$.
\end{prop}
\begin{proof}
It suffices to note that the cokernel of $F-1$ on $R_\sigma$ is generated
freely by the images of $\sigma \cap (\ZZ^n \setminus p\ZZ^n)$.
This yields the first assertion; the second follows by
Theorem~\ref{T:p-limit}.
\end{proof}

\subsection{Projections and sections}

\begin{defn}
A convex cone $\sigma$ is \emph{strictly convex} if for $\bv, \bw \in \sigma$,
$\bv + \bw = 0$ if and only if $\bv = \bw = 0$.
For $\sigma$ a strictly convex cone, let $R'_\sigma$ be the subring of 
$R_\sigma$
consisting of elements $\sum_{\bv} c_\bv [\bv]$ with $c_0 \in \FF_p$.
(Note that strict convexity is needed for this subset to be closed under
multiplication.)
\end{defn}

\begin{prop} \label{P:split2}
Suppose that 
$\sigma$ and $\sigma_0$ are convex cones,
and $\{\sigma_i\}_{i \in I}$ is a (not necessarily finite)
collection of strictly convex cones,
such that $\sigma \setminus \{0\}$ is the disjoint union of
$\sigma_0 \setminus \{0\}$ and the $\sigma_i \setminus \{0\}$.
Then the natural map
\[
\pi_1^p(\Spec R_\sigma) \to
\pi_1^p(\Spec R_{\sigma_0}) 
\times \prod_{i \in I}
\pi_1^p(\Spec R'_{\sigma_i}) 
\]
is an isomorphism.
\end{prop}
\begin{proof}
The argument is as in Proposition~\ref{P:split1}. 
\end{proof}

\begin{defn} \label{D:pi sigma tau}
Let $\sigma, \tau$ be convex cones with $\tau \subseteq \sigma$.
Put $\sigma_0 = \tau$, and choose a 
collection $\{\sigma_i\}_{i \in I}$ of strictly convex cones such that
$\sigma \setminus \{0\}$ is the disjoint union of
$\sigma_0 \setminus \{0\}$ and the $\sigma_i \setminus \{0\}$.
Then the product decomposition given by
Proposition~\ref{P:split2} yields a morphism
\[
\pi_{\sigma,\tau}: \pi_1^p(\Spec R_\tau)
\to \pi_1^p(\Spec R_\sigma)
\]
sectioning the projection $\pi_1^p(\Spec R_\sigma) \to \pi_1^p(\Spec R_\tau)$.
Note that replacing one of the $\sigma_i$ by a disjoint union
does not affect $\pi_{\sigma,\tau}$; in particular, by passing to
a common refinement, we see that this map does not depend at all on
the choice of the $\sigma_i$.
\end{defn}

\begin{prop} \label{P:cut down AS}
Let $\sigma,\tau$ be convex cones with $\tau \subseteq \sigma$.
Let $\rho: \pi_1^p(\Spec R_\sigma) \to \ZZ/p\ZZ$ be the homomorphism
corresponding to the $\ZZ/p\ZZ$-torsor $S = R_\sigma[z]/(z^p - z - x)$ over
$R_\sigma$. Write $x = \sum_{\bv \in \sigma \cap \ZZ^n} c_{\bv} [\bv]$.
Then $\rho \circ \pi_{\sigma,\tau}: \pi_1^p(\Spec R_\tau) \to \ZZ/p\ZZ$
is the homomorphism corresponding to the $\ZZ/p\ZZ$-torsor
\[
R_\tau[z]/(z^p - z - x'), \qquad x' = \sum_{\bv \in \tau \cap \ZZ^n}
c_{\bv}[\bv].
\]
\end{prop}
\begin{proof}
This is an immediate consequence of how
the map in Theorem~\ref{T:p-limit} is constructed.
\end{proof}

We can make the maps $\pi_{\sigma, \tau}$ more explicit in certain cases
as follows.
\begin{defn} \label{D:vlambda}
Let $\lambda: \RR^n \to \RR$ be a linear function. For $\sigma$ a convex
cone, let $v_\lambda$ be the valuation on $R_\sigma$ defined by
\[
v_{\lambda} \left( \sum_{\bv \in \sigma \cap \ZZ^n} c_\bv [\bv] \right)
= \min \{ \lambda(\bv): \bv \in \sigma \cap \ZZ^n, c_\bv \neq 0 \}.
\]
Let $R_{\sigma, \lambda}$ be the completion of $R_\sigma$ with respect 
to $v_\lambda$, and put
\[
\tau_{\sigma,\lambda} = \{\bv \in \sigma: \lambda(\bv) \leq 0\}.
\]
\end{defn}
\begin{prop} \label{P:cutaway}
Set notation as in Definition~\ref{D:vlambda}, and write
$\tau$ for $\tau_{\sigma, \lambda}$.
Then
the natural map $\pi_1^p(\Spec R_{\sigma,\lambda}) 
\to \pi_1^p(\Spec R_{\tau})$ is an isomorphism.
\end{prop}
\begin{proof}
If $z \in R_{\sigma,\lambda}$ and $v_\lambda(z) > 0$, then 
$z + z^p + \cdots$ converges in $R_{\sigma,\lambda}$ to some $y$
satisfying $y^p - y = -z$. Thus the cokernels of $F-1$ on
$R_\tau$ and $R_{\sigma,\lambda}$ are isomorphic, so 
Theorem~\ref{T:p-limit} applies.
\end{proof}

\begin{example}
Simple examples of Proposition~\ref{P:cutaway} are
the facts that $\pi_1^p(\Spec R \llbracket t \rrbracket) \to \pi_1^p(\Spec R)$
and $\pi_1^p(\Spec R((t))) \to \pi_1^p(\Spec R(t^{-1}))$ are isomorphisms.
For a more nontrivial example, take $\sigma$ to be the nonnegative quadrant
in $\RR^2$, and define $\lambda(a,b) = a-b$. Then Proposition~\ref{P:cutaway}
asserts that
\[
\pi_1^p(\Spec (R[xy]\llbracket x \rrbracket)[y])) \to \pi_1^p(\Spec R[xy, y])
\]
is an isomorphism.
\end{example}

\subsection{Heights and representations}

\begin{defn}
A \emph{linear cone} in $\RR^n$ is a convex cone consisting of the
nonnegative scalar multiples of a single nonzero element of $\RR^n$.
For $S$ a set of linear cones and $\sigma$ a convex cone, define
\[
S_\sigma = \{\tau \in S: \tau \subseteq \sigma\}.
\]
\end{defn}

\begin{theorem} \label{T:discrete}
Let $\sigma$ be a convex cone in $\RR^n$, 
let $h$ be a height function on $\calC_{R_\sigma}$ over $R$,
let $\ell$ be a nonzero real number,
and let $G$ be a finite discrete group.
Then there exists a finite set $S$ of linear cones in $\RR^n$,
depending on $\sigma, h, \ell, G$,
such that for any continuous representation
$\rho: \pi_1^p(\Spec R_\sigma) \to G$ with $h(\rho) \leq \ell$
and any convex cone $\tau \subseteq \sigma$, the
image $(\rho \circ \pi_{\sigma,\tau})(\pi_1^p(\Spec R_\tau))$
is determined by $\rho$ and $S_\tau$.
\end{theorem}
\begin{proof}
We first check the case $G = \ZZ/p\ZZ$.
If $\rho$ is trivial, there is nothing to check; otherwise, $\rho$ becomes
trivial upon restriction to $\pi_1^p(\Spec S)$ for some
$\ZZ/p\ZZ$-torsor $S$ over $R_\sigma$. By Proposition~\ref{P:as},
we may write $S = R_{\sigma}[z]/(z^p - z - x)$, and we may choose
$x \in R_{\sigma}$ with support in $\{0\} \cup 
(\sigma \cap (\ZZ^n \setminus p\ZZ^n))$.
Since $h(\rho) \leq \ell$, by (e)
the support of $x$ is contained in a finite set $T$ depending on 
$\sigma,h,\ell$.
By Proposition~\ref{P:cut down AS},
the claim holds if we take $S$ to be the set of linear cones
generated by the elements of $T$; note that this set depends only on
$\sigma, h,\ell$, and not on $\rho$.

We next pass to the general case.
We may assume that $G$ is the image of $\rho$,
so that $G$ is a $p$-group and $\rho$ is surjective; we may also assume
that $G$ is nontrivial.
Let $K$ be the Frattini subgroup of $G$ (the intersection of its maximal
proper subgroups), so that $G/K$ is an elementary abelian $p$-group.
By repeatedly applying the previous paragraph, we obtain a finite set $S_1$
of linear cones, determined by $\sigma, h, \ell$,
such that the image of $\pi_1^p(\Spec R_\tau)$ in $G/K$
is determined by $\rho$ and $(S_1)_\tau$.

We now induct on the size of (the smallest possible choice of) $S_1$.
If $S_1$ is empty, then the image of $\pi_1^p(\Spec R_\tau)$ in $G/K$
is equal to the image of $\pi_1^p(\Spec R_\sigma)$ in $G/K$, namely
$G/K$ itself. Thus the image of $\pi_1^p(\Spec R_\tau)$ in $G$
cannot be contained in any proper subgroup of $G$, and so must equal
$G$; we are thus done in this case.

If $S_1$ is nonempty,
choose a linear cone $T$ in $S_1$; we can then choose strictly convex cones
$\sigma_1, \dots, \sigma_m$ not meeting $T$ such that
$\sigma \setminus T$ is the union of 
$\sigma_1 \setminus \{0\}, \dots, \sigma_m \setminus \{0\}$.
(Namely, draw $n-1$ hyperplanes meeting transversely along $T$,
take the open halfspaces on both sides of each plus one halfspace
containing the negation of $T$, and intersect all of these with
$\sigma$.) We may now apply the induction hypothesis to each
of the $\sigma_i$ (since the analogue of the set $S_1$ has been
reduced by one element) to produce a finite set $S_T$ 
(determined by $\sigma_i, h,\ell,G$) such that if 
$T \not\subseteq \tau$, then $\rho(\pi_1^p(\Spec R_\tau))$ is determined
by $\rho$ and $(S_T)_\tau$. Let $S$ be the union of the $S_T$; this
has the desired property because if every $T$ lies in $\tau$,
then the image of $\pi_1^p(\Spec R_\tau)$ in $G/K$ must equal $G/K$,
so as in the base case, $\rho(\pi_1^p(\Spec R_\tau)) = G$.
\end{proof}
\begin{cor}
Let $\sigma$ be a convex cone in $\RR^n$, 
let $h$ be a height function on $\calC_{R_\sigma}$ over $R$,
let $\ell$ be a nonzero real number,
and let $G$ be a finite discrete group.
Then there exists a finite set $S$ of linear cones in $\RR^n$,
depending on $\sigma, h, \ell, G$,
such that for any continuous representations
$\rho_1, \rho_2: \pi_1^p(\Spec R_\sigma) \to G$ with $h(\rho_1), h(\rho_2) 
\leq \ell$ and any convex cone $\tau \subseteq \sigma$, whether or not
the restrictions of $\rho_1$ and $\rho_2$ to $\pi_1^p(\Spec R_\tau)$
are isomorphic
is determined by $\rho_1, \rho_2, S_\tau$.
\end{cor}
\begin{proof}
Embed $G$ into a linear group over a field of characteristic zero,
and apply Theorem~\ref{T:discrete} to the representation
$\rho_1^{\dual} \times \rho_2: \pi_1^p(\Spec R_\sigma) \to G^T \times G$.
(Here $\rho_1^{\dual}$ denotes the contragredient representation and $G^T$
denotes $G$ with its linear embedding transposed.)
\end{proof}

The next corollary is sufficiently useful in its own right that we
have promoted it to a theorem.
\begin{theorem} \label{T:splits2}
Fix a convex cone $\sigma$.
For $T = \{\tau_1, \dots, \tau_m\}$ a collection of distinct linear cones
contained in $\sigma$,
let $G_T$ be the limit of the diagram consisting of the arrows
$\pi^p_1(\Spec R_{\tau_i}) \to \pi_1^p(\Spec R)$ for $i=1, \dots, m$.
View the $G_T$ as an inverse system via the natural maps $G_{T'} \to G_T$
for $T \subseteq T'$.
Then $\pi_1^p(\Spec R_\sigma)$ is the inverse
limit of the $G_T$.
\end{theorem}
A weaker but coordinate-free variant of the Theorem~\ref{T:splits2}
is the following.
\begin{cor}
Fix a convex cone $\sigma$.
For $S = \{R_1, \dots, R_m\}$ a set whose elements are subalgebras
of $R_{\sigma}$ of transcendence degree $1$ over $R$, let
$G_S$ be the limit of the diagram consisting of the arrows
$\pi^p_1(\Spec R_i) \to \pi_1^p(\Spec (R_i \cap R_j))$ for $i,j=1, \dots, m$.
View the $G_S$ as an inverse system via the natural maps $G_{S'} \to G_S$
for $S \subseteq S'$.
Then $\pi_1^p(\Spec R_\sigma)$ is the inverse limit of the $G_S$.
\end{cor}
Finally, it is worth saying in simple terms what Theorem~\ref{T:splits2}
says about affine spaces.
\begin{cor}
For $n$ a positive integer $n$, take $x_1, \dots, x_n$ to be coordinates on
$\AAA^n_R$. Then the group $\pi_1^p(\AAA^n_R)$ is the inverse limit
of the groups $\pi_1^p(\Spec R[x_1^{a_1} \cdots x_n^{a_n}])$
over all coprime $n$-tuples $(a_1, \dots, a_n)$ of nonnegative integers.
\end{cor}

\section{Complements on height functions}

To conclude, we point out that the somewhat mysterious height functions
that we have been using can be made quite explicit on
affine toric varieties. The main result here
is Theorem~\ref{T:height splits}, which gives a relatively simple formula
for a function which can be verified (Corollary~\ref{C:strong height})
to be a height function. 

Note that we use Theorem~\ref{T:splits2}
in the course of the proof; we do not know whether it is possible
to prove Theorem~\ref{T:height splits} directly, then 
short-circuit the proof of Theorem~\ref{T:discrete} and its consequences
around the discussion of general height functions. 
Doing so might necessitate establishing a relationship between
ramification theory for local fields with
imperfect residue field (see Remark~\ref{R:imperfect}); such a relationship
might have the effect of clarifying the ramification theory in some cases.

\subsection{Some explicit height functions}

In the situation we have been considering, we can write down some height
functions explicitly.
\begin{convention}
Throughout this section, let $R = k$ be an algebraically closed $p$-field.
\end{convention}

\begin{defn}
For $\lambda: \RR^n \to \RR$ a nonzero linear function defined over
$\QQ$ (i.e., it carries $\QQ^n$ to $\QQ$), 
let $m_\lambda$ be the unique rational number such that
$m_\lambda \lambda(\ZZ^n) = \ZZ$,
let $H_\lambda$ denote the hyperplane $\{\bv \in \RR^n: \lambda(\bv) = 0\}$,
and let $K_\lambda$ denote the perfection of the fraction field of
$R_{H_\lambda}$.
Let $\widehat{R}_\lambda$ denote the completion of
$R_{\RR^n}$ with respect to $v_{-\lambda}$,
and put
\[
Q_{\lambda} = \widehat{R}_\lambda \otimes_{R_{H_\lambda}} K_\lambda;
\]
we may then view $Q_\lambda$ as a power series field in one variable
over the perfect field $K_\lambda$, with valuation
$m_\lambda v_{-\lambda}$.
Given a $p$-typical extension $S$ of $R_{\RR^n}$, define
\[
c_\lambda(S) = \frac{1}{m_\lambda} 
b((S \otimes_{R_{\RR^n}} Q_{\lambda})/Q_{\lambda}),
\]
where $b$ denotes the highest break function
(of Definition~\ref{D:highest break}).
\end{defn}

As in Remark~\ref{R:onedim}, $b_\lambda$ is not a height function for
$p$-typical extensions of $R_{H_\lambda}$. However, we can use the
functions $b_\lambda$ to construct height functions on smaller cones.

\begin{defn}
Given a convex cone $\sigma$, define the \emph{dual cone} $\sigma^\dual 
\subseteq (\RR^n)^\dual$ as the set of linear functions
$\lambda: \RR^n \to \RR$ such that $\lambda(\bv) \geq 0$ for all $\bv 
\in \sigma$. We say $\sigma$ is \emph{very convex} if 
$\sigma^\dual$ has nonempty topological interior relative to $(\RR^n)^\dual$;
if $\sigma$ is very convex, then it is strictly convex.
\end{defn}

\begin{defn} \label{D:height}
Let $\sigma$ be a nontrivial very convex cone,
and let $U \subseteq \sigma^\dual \setminus \{0\}$ be a subset
open in $(\RR^n)^\dual$.
Define the function $h_{U}$ on $\calC_{R_\sigma}$ by
\[
h_{U}(S) = \sup_{\lambda \in U \cap (\QQ^n)^\dual} 
\{c_{\lambda}(S \otimes R_{\RR^n})\}.
\]
For $\lambda$ in the interior of $\sigma^\dual \setminus \{0\}$, define
\[
h_\lambda(S) = \limsup_U h_U(S),
\]
the limit being taken over the direct system of open neighborhoods of $\lambda$
in $\sigma^\dual \setminus \{0\}$.
For $\rho: \pi_1^p(R_\sigma) \to G$ a continuous representation to
a discrete group, put $h_U(\rho) = h_U(S)$ and $h_\lambda(\rho)
= h_\lambda(S)$, for $S \in \calC_{R_\sigma}$ connected and
chosen so that $\ker(\rho) = \pi_1^p(S)$.
\end{defn}

We first work out how $h_U$ works on linear cones. First, we bundle together
some hypotheses.
\begin{hypothesis} \label{H:linear}
Let $\sigma \subseteq \RR^n$ be a linear cone with $\ZZ^n 
\cap \sigma \neq \{0\}$, and put $\tau = -\sigma \cup \sigma$.
Let $\widehat{R}_\sigma$ be the completion of $R_\tau$ with respect to
$v_{-\lambda}$, for some nonzero linear functional $\lambda: \RR^n \to \RR$
defined over $\QQ$ which is positive on $\sigma \setminus \{0\}$.
Let $\rho: \pi_1^p(\Spec R_{\sigma}) \to G$ be a continuous
representation to a discrete group.
Note that $\widehat{R}_\sigma$ is a power series field over $k$,
and that it depends only on $\sigma$, not on $\lambda$; we may thus sensibly
speak of the highest break $b(\rho)$.
\end{hypothesis}

\begin{lemma}\label{L:specialize1}
Under Hypothesis~\ref{H:linear},
let $\lambda: \RR^n \to \RR$ be a nonzero linear function defined over
$\QQ$, such that $\lambda$ is positive on $\sigma \setminus \{0\}$.
Put $d = [\ZZ^n: (\ZZ^n \cap H_\lambda) \times (\ZZ^n \cap \tau)]$,
and let $d'$ be the prime-to-$p$ part of $d$.
Then
\[
c_{\lambda}(\rho) = \frac{d'}{m_\lambda} b(\rho).
\]
\end{lemma}
\begin{proof}
We first note that the desired equality holds when
$d = 1$: it is the comparison between the highest
break of a representation of $\pi_1^p$ of a power series ring over a field,
and the same representation after pulling back by extending the constant field.

We next note that if we repeat the construction of $c_\lambda(\rho)$ 
with $\ZZ^n$
replaced by the larger lattice 
$(\ZZ^n \cap H_\lambda) \times \frac{1}{d}(\ZZ^n \cap \tau)$,
then $m_\lambda$ and $c_\lambda(\rho)$ remain unchanged.
However, by Definition~\ref{D:highest break},
$b(\rho)$ gets multiplied by $d'$. Now appealing to the $d=1$ case
yields the desired result.
\end{proof}

\begin{cor} \label{C:linear union}
With notation as in Definition~\ref{D:height} and
Hypothesis~\ref{H:linear},
let $\bv$ be the smallest nonzero element of $\ZZ^n \cap \sigma$.
Then
\[
h_U(\rho) = b(\rho)  \sup_{\lambda \in U} \{ \lambda(\bv) \},
\qquad
h_\lambda(\rho) = b(\rho) \lambda(\bv).
\]
In particular,
\[
h_U(\rho) = \sup_{\lambda \in U} \{ h_\lambda(\rho)\}.
\]
\end{cor}
\begin{proof}
With notation as in 
Lemma~\ref{L:specialize1}, note that
\[
d = [\lambda(\ZZ^n):\lambda(\ZZ^n \cap \sigma)]
= m_\lambda \lambda(\bv).
\]
By Lemma~\ref{L:specialize1}, we then have
\[
c_\lambda(\rho) \leq b(\rho) \lambda(\bv),
\]
with equality for any $\lambda$ for which $d$ is not divisible by $p$.
Such $\lambda$ are dense in any $U$, so the desired results follow.
\end{proof}

We now treat general cones.
\begin{defn}
For $\sigma$ a very convex cone, $\tau$ a convex cone contained in $\sigma$,
and $\rho: \pi_1^p(\Spec R_\sigma) \to G$ a continuous representation
to a discrete group,
let $\rho_\tau$ be the pullback of $\rho$ along the maps
\[
\pi_1^p(\Spec R_\sigma) \to \pi_1^p(\Spec R_\tau) \to \pi_1^p(\Spec R_\sigma),
\]
where the second map is as in Definition~\ref{D:pi sigma tau}.
\end{defn}

\begin{lemma}\label{L:specialize}
With notation as in Definition~\ref{D:height}, 
let $T$ be a linear cone contained in $\sigma$ such that
\[
d = [\ZZ^n:(\ZZ^n \cap H_\lambda) \times (\ZZ^n \cap (T \cup -T))]
\]
is coprime to $p$.
Then
\[
c_{\lambda}(\rho) \geq c_{\lambda}(\rho_T).
\]
\end{lemma}
\begin{proof}
As in the proof of Lemma~\ref{L:specialize1}, we may reduce to the case $d=1$;
then Proposition~\ref{P:Laumon} yields the claim.
\end{proof}
\begin{cor}
With notation as in Lemma~\ref{L:specialize}, we have
\[
h_U(\rho) \geq h_U(\rho_T),
\qquad
h_\lambda(\rho) \geq h_\lambda(\rho_T).
\]
\end{cor}
\begin{proof}
By Lemma~\ref{L:specialize1}, 
we may compute $h_U(\rho_T)$ by taking the supremum defining it
only over $\lambda$ as in Lemma~\ref{L:specialize} (i.e., the
$\lambda$ for which $d=d'$ in Lemma~\ref{L:specialize1}).
Then Lemma~\ref{L:specialize} yields the first desired inequality;
the second follows by taking limits.
\end{proof}

We now have the following fairly explicit description of the functions
$h_U$ and $h_\lambda$, in terms of linear cones.
\begin{theorem} \label{T:height splits}
With notation as in Definition~\ref{D:height}, we have
\begin{equation} \label{eq:height splits}
h_U(\rho) = \sup_T \{ h_U(\rho_T)\},
\qquad
h_{\lambda}(\rho) = \sup_T \{ h_{\lambda}(\rho_T) \},
\end{equation}
the suprema taken over all linear cones $T \subseteq \sigma$.
\end{theorem}
\begin{proof}
In each case, the left side is greater than or equal to the right by 
Lemma~\ref{L:specialize}. Conversely,
by Theorem~\ref{T:splits2}, we can present $\rho$ inside the tensor
product of the $\rho_T$ over finitely many $T$, and so the right
side is greater than or equal to the left.
\end{proof}

\subsection{More on the explicit height functions}

Theorem~\ref{T:height splits} makes it easy to verify many natural properties
of the $h_\lambda$, including the fact that they are actually height functions.
We present these as a series of corollaries.
\begin{convention}
Throughout this section, retain notation as in Definition~\ref{D:height}.
\end{convention}

\begin{cor}
We have
\[
h_U(\rho) = \sup_{\lambda \in U} \{ h_\lambda(\rho)\}.
\]
\end{cor}
\begin{proof}
Applying Theorem~\ref{T:height splits} and Corollary~\ref{C:linear union},
we have
\[
h_U(\rho) = \sup_{T \subseteq \sigma} \{h_U(\rho_T)\}
= \sup_{T \subseteq \sigma, \lambda \in U} \{h_\lambda(\rho_T)\}
= \sup_{\lambda \in U} \{h_\lambda(\rho)\}.
\]
\end{proof}

\begin{cor}
If $\lambda \in \sigma^\dual \setminus \{0\}$ is defined over $\QQ$,
then $h_\lambda(\rho) \in \QQ$.
\end{cor}
\begin{proof}
Apply Theorem~\ref{T:height splits} and note that only finitely many
terms in the supremum in \eqref{eq:height splits}
are nonzero thanks to Theorem~\ref{T:splits2}.
Then apply Corollary~\ref{C:linear union} to deduce that each nonzero
term in the supremum is rational.
\end{proof}

\begin{cor}
Suppose that $\tau$ is a convex cone with $\tau \subseteq \sigma$. Then
\[
h_\lambda(\rho) \geq h_\lambda(\rho_\tau).
\]
\end{cor}
\begin{proof}
Apply Theorem~\ref{T:height splits}, and note that the supremum
defining $h_\lambda(\rho_\tau)$ is simply the same supremum 
as in \eqref{eq:height splits}, but restricted to $T \subseteq \tau$.
\end{proof}
\begin{cor}
Suppose that $\lambda$ and $\kappa$ both
lie in the interior of $\sigma^\dual \setminus \{0\}$,
and that $\lambda(\bv) \geq \kappa(\bv)$ for all $\bv \in \sigma$.
Then
\[
h_\lambda(\rho) \geq h_\kappa(\bv).
\]
\end{cor}
\begin{proof}
By Theorem~\ref{T:height splits}, it suffices to check this for
$\sigma$ a linear cone, in which case it follows from
Corollary~\ref{C:linear union}.
\end{proof}

\begin{cor} \label{C:height of AS}
Suppose that 
$S = R_\sigma[z]/(z^p - z - x)$, where the support $V$ of
$x$ is contained in $\sigma \cap (\ZZ^n \setminus p\ZZ^n)$.
Then
\[
h_{\lambda}(S) = \sup_{\bw \in V} \{ v_\lambda(\bw)\}.
\]
\end{cor}
\begin{proof}
Apply Theorem~\ref{T:height splits} to reduce to the case
where $\sigma$ is linear. Then apply
Corollary~\ref{C:linear union} and Lemma~\ref{L:as1}.
\end{proof}
\begin{cor} \label{C:strong height}
Each of the functions
$h_\lambda$ and $h_U$ is a strong height function
on $\calC_{R_\sigma}$ over $R = k$.
\end{cor}
\begin{proof}
Conditions (a)-(d) are straightforward, while condition (e) for $d=1$
follows from Corollary~\ref{C:height of AS}; the claim then follows
by Proposition~\ref{P:height1}.
\end{proof}

\end{document}